\def\IR{{\mathbb R}}

\def\IC{\mathbb C} 
\def\ID{{\mathbb D}}
\def\oC{\hat{\IC}}

\def\zbar{{\overline{z}}}

\documentclass[12pt]{article} 

\usepackage[psamsfonts]{amssymb} 
\usepackage{amsfonts,amsmath}

\usepackage{epsfig,multicol} 
\DeclareMathOperator{\loc}{loc}

\newtheorem{theorem}{Theorem}
\newtheorem{lemma}{Lemma}

\title{Quasiregular families bounded in $L^p$ and elliptic estimates.}
\author{Aimo Hinkkanen and Gaven Martin}  
\begin{document}

\maketitle

\begin{abstract}  
We prove that a family ${\cal F}$ of quasiregular mappings of a domain $\Omega$ which are uniformly bounded in $L^p$ for some $p>0$ form a normal family.  From this we show how an elliptic estimate on a functional differences implies all directional derivatives, and thus the complex gradient to be quasiregular.  Consequently the function enjoys much higher regularity than apriori assumptions suggest.
\end{abstract}

\section{Introduction}    

The governing equations of geometric function theory and the theory of quasiconformal mappings,  Teichm\"{u}ller spaces and so forth are the Beltrami equations and their nonlinear counterparts,  see for instance  \cite{AIM,IM,G1,G2,H} and elsewhere.  Beltrami equations come in several different flavours.  As examples,  let $\Omega\subset \IC$ be a domain and let $f:\Omega\to\IC$ be a mapping of Sobolev class $W^{1,1}_{loc}(\Omega)$ consisting of functions whose first derivatives are locally integrable;  then we have
\begin{itemize}
\item $\IC$-linear:  $f_\zbar=\mu(z) f_z$,  with ellipticity estimate $\|\mu\|_{L^\infty(\Omega)}<1$;
\item $\IR$-linear:  $f_\zbar=\mu(z) f_z+\nu(z) \overline{f_z}$,  with ellipticity estimate \[ \|\,|\mu|+|\nu|\,\|_{L^\infty(\Omega)}<1;\]
\item Autonomous: $f_\zbar={\cal A}(f_z)$,  with ellipticity estimate: there is $k<1$ so that for all $\zeta,\eta\in\IC$ \[|{\cal A}(\zeta)-{\cal A}(\eta)|\leq k |\zeta-\eta|  ;  \] 
\item Fully nonlinear: $f_\zbar={\cal H}(z,f,f_z)$,  with ellipticity estimate: there is $k<1$ so that for all $z\in \Omega$ and for  all $w,\zeta,\eta\in\IC$, we have
 \[|{\cal H}(z,w,\zeta)-{\cal H}(z,w,\eta)|\leq k |\zeta-\eta|  ,  \] 
with additional conditions on ${\cal H}$,  see \cite[Chapters 7 \& 8]{AIM}.
\end{itemize}

Each of these equations has a seminal application and they are all inter-related.  The apriori assumption that $f\in W^{1,1}_{loc}(\Omega)$ is so that we can even speak of $f$ as a ``solution''.  Without stronger assumptions on $\mu$ or ${\cal H}$ not much can be said,  but note for instance that $\mu=0$ on an open set implies $f$ is holomorphic on that set --- Weyl's Lemma.  The higher regularity theory of these equations typically assumes more on $f$,  for instance $f\in W^{1,q}_{loc}(\Omega)$ for some $q$ with $1<q \leq 2 $ usually depending on the ellipticity constant $k$,  and in return delivers a far nicer outcome,  $f\in W^{1,p}_{loc}(\Omega)$ for some $p > 2$,  again depending on $k$.  Astala's theorem \cite{Ast} gives the optimal result in the $\IC$-linear case and can be used to analyse other cases.  Questions of existence and uniqueness are fairly well understood through the topological properties of these mappings and the well known Sto\"{\i}low factorisation theorem, \cite{St} and the references therein,  see also \cite[\S 5.5 \& \S 6.1]{AIM}.  However there are intriguing subtleties in the nonlinear case,  see \cite{AC1,AC2}.

\medskip

In this paper we seek general methods to go beyond the $W^{1,p}_{loc}$ regularity to seek $W^{2,p}_{loc}$ estimates,  see our Theorem \ref{thm2}.  Such estimates have been found before in special cases,  for instance in  the study of the autonomous equations (e.g.,  \cite{AC1}),  and these estimates have important applications (e.g., \cite{Mar}) and serve as a bootstrap for $C^\infty$-regularity.  It is noteworthy that an elliptic estimate such as  (\ref{1}) below implies that the derivative $f_z$  is quasiregular (a solution to an elliptic Beltrami equation).

\section{Main results}

We denote the set of real numbers by ${\mathbb R}$, the set of complex numbers by ${\mathbb C}$, and the unit disk by $ {\mathbb D} = \{ z\in {\mathbb C} \colon |z|<1\}$. We write ${\mathbb D}(z,r)=\{w\in {\mathbb C} \colon |w-z|<r\}$ for the open disk with centre $z\in {\mathbb C}$ and radius $r>0$. Thus $ {\mathbb D} ={\mathbb D}(0,1)$. We denote the Riemann sphere by $\overline{ {\mathbb C} } = {\mathbb C} \cup \{\infty\}$. 

The main results that we present here are the following.
\begin{theorem}\label{Lpbounded}
  Let $\Omega$ be a domain in the complex plane ${\mathbb C}$. Let $p$ be a real number with $p > 0$. Let ${\mathcal F}$ be a family of $K$--quasiregular mappings $f:\Omega\to\IC$ which is uniformly bounded in $L^p(\Omega)$.   Then the family ${\mathcal F}$ is precompact,  every sequence contains a locally uniformly convergent subsequence and each limit function is $K$-quasiregular (possibly constant) with values in $\IC$.
\end{theorem}

We use this theorem to prove the following result which establishes very strong regularity from a standard elliptic type estimate.

\begin{theorem}\label{thm2}  Suppose that  $0\leq k < 1$ and that $f:\Omega \to \IC$ is a function in the  Sobolev class $W^{1,p}_{\loc}(\Omega)$ for some $p>1+k$. Let 
$a:\Omega\to \IR_+$ be a continuous function such that  $0 < a(z)\leq {\rm dist}(z,\partial \Omega)$ for all $z\in\Omega$.  
Suppose further that for a.e.\ $z\in \Omega$ and for every $\zeta$ with $|\zeta|=1$, the function $f$ satisfies the elliptic estimate 
\begin{equation}\label{1}
|f_\zbar(z+t\zeta)-f_\zbar(z) | \leq k \,  |f_z(z+t\zeta)-f_z(z) | 
\end{equation}
for all $t$ with $0< t< a(z)$.  
Then the following hold:
\begin{enumerate}
\item $f\in W^{2,q}_{loc}(\Omega)$ for all $q<1+1/k$.  
\item Each member of the $\IR$-linear family
\[ \{ a f_x(z)+b f_y(z) : a,b\in \IR \} \]
is a $\frac{1+k}{1-k}$-quasiregular mapping, possibly constant.
\item There are measurable functions $\mu,\nu:\Omega\to \IC$ with $|\mu|+|\nu|\leq k$ a.e.\ in $\Omega$ such that both directional derivatives $f_x$ and $f_y$ satisfy the $\IR$-linear Beltrami equation,
\[ h_\zbar = \mu(z) h_z +\nu(z) \overline{h_z}, \hskip20pt h\in \{f_x,f_y\}  .  \]
\item The complex $z-$derivative $f_z$ is itself quasiregular and satisfies, with $\mu$ and $\nu$ as in part 3 above, the $\IR-$linear equation
\begin{equation} \label{rrr}
 h_\zbar = \frac{\mu(z)}{1-|\nu(z)|^2}  h_z + \frac{\overline{\mu(z)}\nu(z)}{1-|\nu(z)|^2}  \overline{h_z}  \hskip20pt \hbox{ where } h = f_z ,  
\end{equation}
and thus $f_z \in W^{1,q}_{loc}(\Omega)$ for all 
\[ q <  1+1/k',  \hskip15pt  \hbox{ where }  k'= \left\| \frac{|\mu|}{1-|\nu|} \right\|_{L^{\infty}(\Omega)} \leq k .  \]
\end{enumerate}
\end{theorem}

The next result concerns the tangent cone of a quasiregular mapping and H\"{o}lder regularity. We denote the Lebesgue area measure by $m$ and integration with respect to the complex variable $z$ using the measure $m$ by $dm(z)$.

\begin{theorem}\label{thm3}  Let $f:\ID\to\IC$ be quasiregular with $f(0)=0$.  Suppose that for some $p,q>0$ and for all $\epsilon$ sufficiently small we have
\begin{equation} \label{integral}
 \int_{\ID(0,\epsilon)} |f(z)|^p \, dm(z) \leq C \epsilon^{2+q} 
 \end{equation}
where $C$ is an absolute constant.  Then the family of quasiregular maps 
\[ {\cal F}=\left\{ \frac{1}{\lambda^{q/p}} f(\lambda z) : \lambda \in \ID\setminus \{0\} \right\} \]
is precompact. Every sequence from ${\mathcal F}$  contains a subsequence converging locally uniformly in $\ID$ to a quasiregular mapping,  or to a constant (possibly $\infty$). 
\end{theorem}
The Sto\"{\i}low factorisation theorem \cite[Theorem 5.5.1, p.~179]{AIM} together with the H\"{o}lder continuity properties of $K$-quasiconformal mappings tells us,  since $f(0)=0$,  that if the local index  of $f$ at $0$ is $n\geq 1$,  then we have the {\em a priori} bound
\[ \int_{\ID(0,\epsilon)} |f(z)|^p \, dm(z) \leq  C \int_{\ID(0,\epsilon)} |z|^{np/K} \, dm(z)  = \frac{2\pi C}{2+np/K} \;  \epsilon^{2+np/K}. \]
So no matter which exponent $p>0$ is chosen,  there is always an exponent $q>0$ such that (\ref{integral}) holds.

\section{Proof of Theorem \ref{thm2}} 

Let the assumptions of Theorem \ref{thm2} be satisfied. In particular, $p$ is a fixed exponent with $p>1+k$ of which we assume that $f:\Omega \to \IC$ lies in $W^{1,p}_{\loc}(\Omega)$.

Let $U$ and $V$ be relatively compact subdomains of $\Omega$ with $\overline{U}\subset V$. Then $a(z)\geq \varepsilon=\varepsilon(V)>0$ on $V$. We furthermore choose $\varepsilon$ so that $ 2 \varepsilon < {\rm dist}\, (V, \partial \Omega)$. Write $ V' = \{ z\in {\mathbb C} \colon {\rm dist}\, (z,V)<\varepsilon\}$ so that $\overline{V} \subset V'$ and $\overline{V'}  \subset \Omega$. 

Fix $t$ and $\zeta$ so that $0<t<\varepsilon$ and $|\zeta|=1$. Then $g(z)= f(z+t\zeta)-f(z)$ is well defined in $V'$ and $g\in W^{1,p}(V')$. 

Due to this and (\ref{1}),  \cite[Theorem 5.4.2, p.~175]{AIM} implies that for each $q\in (1+k,1+1/k)$, and hence for $q=2$, we have $g\in W^{1,q}_{\loc}(V')$. In particular, $g$ is continuous and $g$ is $K-$quasiregular in $V'$, where $K=(1+k)/(1-k)$.

We denote the $2\times 2$ complex derivative matrix of $f$ by $Df=Df(z)$ and its matrix (operator) norm by $|Df|$.  If $\zeta\in {\mathbb C}$ is viewed as a column vector with two real components, we write $(Df)\cdot \zeta$ for the column vector that arises as the product of $Df$ and $\zeta$. Since $f$ is absolutely continuous on lines, for almost every $z\in V$, we have
$$
g(z) = f(z+t\zeta)-f(z) =\int_0^{t} Df(z+u\zeta) \cdot \zeta \, du
$$
so that
$$
| g(z) | \leq \int_0^{t} |  Df(z+u\zeta) | \, du  
$$
and hence by H\"{o}lder's inequality, for every $q\in (1+k,1+1/k)$, we have 
$$
| g(z) |^q  \leq t^{q-1} \int_0^{t} |  Df(z+u\zeta) |^q \, du  .
$$
It follows that with the notation 
\[ \varphi_t(z) = \frac{1}{t} \big[ f(z+t\zeta)-f(z) \big]  ,   \]
we have
\begin{eqnarray} \label{qq}
\int_V | \varphi_t(z)  |^q \, dm(z)  & \leq &
\frac{1}{t} \int_V \int_0^{t} |  Df(z+u\zeta) |^q \, du  \, dm(z) 
\notag
\\
& \leq &  
\int_{V'} |Df(z) |^q \, dm(z) . 
\end{eqnarray} 
The importance of this estimate is that the right hand side does not depend on~$t$.

The mapping $\varphi_t$ is also $K-$quasiregular in $V'$. Fix $q\in (1+k,1+1/k)$. 
The inequality (\ref{qq}) implies that the family of $K-$quasiregular mappings $\{ \varphi_t \colon 0<t<\varepsilon, \, |\zeta|=1 \}$ is uniformly bounded in $L^q(V)$.  Theorem \ref{Lpbounded}  establishes the precompactness of this family.  In particular, for each fixed $\zeta$ and each sequence $t_j\to 0$, there is a subsequence $t_{j_k}$ such that the limit  $\lim_{k\to\infty} \varphi_{t_{j_k}}$ exists and is $K$-quasiregular (possibly constant). Since $f$ is differentiable almost everywhere, such a limit is equal almost everywhere to the directional derivative
$\partial_{\zeta}f$  of $f$ in the direction $\zeta$.  Thus all directional derivatives of $f$ are $K-$quasiregular. The directional derivatives of $f$ form an $\IR$-linear family of  mappings,  
\begin{equation} \label{real}
 \partial_{\zeta} f = a f_x + b f_y, \hskip15pt \zeta = a+ib,  \hskip5pt a^2+b^2=1  .  
 \end{equation}
Thus, in $V$, the family of all real multiples of directional derivatives (for all real $a$ and $b$, without the restriction $a^2+b^2=1$ in (\ref{real})) of $f$ forms an $\IR$-linear family of quasiregular maps. Since $V$ is arbitrary, the same conclusion holds in $\Omega$. This proves part 2 of Theorem \ref{thm2}. 

The same argument using the function $\varphi_h =   \frac{1}{h} \big[ f(z+h)-f(z)\big]$, where $0<|h|<\varepsilon$, shows that for every sequence $h_j\to 0$ there is a convergent subsequence 
$\varphi_{h_{j_k}} $ whose limit is $K-$quasiregular, possibly constant. The set of all such limits forms, by definition, the tangent cone to $f$, and therefore consists entirely of $K-$quasiregular mappings.

We next appeal to the Caccioppoli type estimate \cite[Theorem 5.4.2, p.~175]{AIM}, choosing the Lipschitz function $\eta$ in that theorem to have compact support in $\Omega$ and to  satisfy $\eta\equiv 1$ on $U$ and $\eta=0$ outside $V$. This tells us that for every $q\in (1+k,1+1/k)$, we have 
\begin{equation} \label{cacc}
 \int_U |D \varphi_t(z) |^q \, dm(z)  \leq C_{U,V}(q)  \int_V |\varphi_t(z)|^q \, dm(z)  ,  \end{equation} 
where the positive constant $C_{U,V}(q)$ depends only on $q$, and on $U$ and $V$ through the specific choice of $\eta$. 

At every point $w\in V'$ where $\varphi_t$ is differentiable, and hence almost everywhere in $V'$, we have 
$$
| (\varphi_t)_x(w)| \leq |D \varphi_t(w) |,\quad | (\varphi_t)_y(w)| \leq |D \varphi_t(w) | ,
$$
so that
$$
| (\varphi_t)_z(w)|  =  \left|  \frac{1}{2}  \left(  \frac{ \partial  \varphi_t } { \partial x   }   -  i  \frac{ \partial  \varphi_t } { \partial y   }  \right)      \right|       \leq  |D \varphi_t(w) | 
$$
while by (\ref{1}), we have
$$
| (\varphi_t)_{\overline{z}}(w)| \leq  k | (\varphi_t)_z(w)| \leq  k  |D \varphi_t(w) | .
$$
Hence, also by (\ref{qq}) and (\ref{cacc}),  for every $q\in (1+k,1+1/k)$, we have 
\begin{eqnarray} \label{6}
&{}& 
\int_U \left|  \frac{f_{\zbar}(z+t\zeta)-f_{\zbar}(z)}{t}  \right|^q \, dm(z)  + \int_U \left|  \frac{f_z(z+t\zeta)-f_z(z)}{t} \right|^q   \, dm(z)   
\notag \\
& \leq & 2 C_{U,V}(q)  \int_V   \left| \frac{f(z+t\zeta)-f(z)}{t}  \right|^q  \, dm(z)
  \notag \\
   &\leq & 
     2 C_{U,V}(q)  \int _{V'}  |Df(z)|^q \, dm(z)  .  
\end{eqnarray}

Since $f_x$ and $f_y$ are $K-$quasiregular in $\Omega$, as we saw earlier, it follows that  for each $q\in (1+k,1+1/k)$, the functions $f_x$ and $f_y$ belong to the Sobolev space $W^{1,q}_{loc}(\Omega)$, again by \cite[Theorem 5.4.2, p.~175]{AIM}. Therefore also $f_z$ and $f_\zbar$ belong to  $W^{1,q}_{loc}(\Omega)$. Hence the directional derivatives of $f_z$ and of $f_\zbar$ exist almost everywhere in $\Omega$. 

Fix $q\in (1+k,1+1/k)$ and fix $\zeta$ with $|\zeta|=1$, and for $t$ with $0<t<\varepsilon$, define  
$$F_t(z)= (f_z(z+t\zeta)-f_z(z))/t$$ 
for $z\in V'$. Then $F_t\in W^{1,q}(V')$. By (\ref{6}), 
$$ \int_U | F_t |^q \leq C \int_{V'}  |Df|^q $$ for all these $t$, where $C=2 C_{U,V}(q)$. We have $\lim_{t\to 0} F_t(z) = \partial_{\zeta}f_z(z)$ for almost every $z\in U$. By Fatou's lemma, we have $\partial_{\zeta}f_z \in L^q(U)$. Similarly, $\partial_{\zeta}f_{\overline{z}} \in L^q(U)$. Taking $\zeta=1$ and $\zeta=i$ and then taking linear combinations of the resulting directional derivatives, we find that also each of the functions $f_{zz}$, $(f_z)_{ \overline{z} }$,  $(f_{ \overline{z}})_{z}$, and $(f_{ \overline{z}})_{ \overline{z} }$ belongs to $L^q(U)$, and there is a uniform upper bound for their norms in $L^q(U)$. Since $U$ is arbitrary and $q$ is only subject to $q\in (1+k,1+1/k)$, it follows that $f\in W^{2,q}_{loc}(\Omega)$ for each $q\in (0,1+1/k)$. This proves part 1 of Theorem \ref{thm2}.

Associated to an $\IR$-linear family $ {\cal F}$ of $\frac{1+k}{1-k}-$quasiregular mappings of a domain $\Omega$,  there are measurable functions $\mu,\nu:\Omega\to \IC$ such that each $h\in {\cal F}$ satisfies the $\IR$-linear equation
\begin{equation} h_\zbar = \mu(z) h_z + \nu(z) \overline{h_z}, \hskip20pt \hbox{ a.e.\ in }\;\; \Omega \label{rlinear} \end{equation}
with the elliptic bounds $|\mu(z)|+|\nu(z)|\leq k$ for almost every $z\in \Omega$. For this result, see \cite{Boj}.  In \cite{JJ} there is a further discussion of this fact, and the uniqueness of $\mu$ and $\nu$ is also proved.  Therefore  $f_x$ and $f_y$ satisfy the equation (\ref{rlinear}).  This proves part 3 of Theorem \ref{thm2}.

Applying (\ref{rlinear}) with $h=f_x$ and $h=f_y$ as appropriate, we find that
\begin{eqnarray*}
2 (f_z)_\zbar & = & (f_x - i f_y)_\zbar  = \mu(z) (f_x)_z+\nu(z) \overline{(f_x)_z} - i\big[ \mu(z) (f_y)_z+\nu(z) \overline{(f_y)_z} \big]\\
& = &   \mu(z) (f_z)_x +\nu(z) \overline{(f_z)_x} - i \big[ \mu(z) (f_z)_y +\nu(z) \overline{(f_z)_y} \big] \\
& = &  2  \mu(z) (f_z)_z + 2 \nu(z) \overline{(f_z)_\zbar}  
\\
& = &  2 \mu(z) (f_z)_z + 2 \nu(z)  \overline{\big[\mu(z) (f_z)_z +\nu(z) \overline{(f_z)_\zbar} \big] }   
\end{eqnarray*}
Hence we obtain the following equation for $f_z$:
\begin{eqnarray*}
(1-|\nu(z)|^2) (f_z)_\zbar  & = &   \mu(z) (f_z)_z +\nu(z)  \overline{\mu(z)}\; \overline{(f_z)_z }   
\end{eqnarray*}
This proves (\ref{rrr}) in part 4 of Theorem~\ref{thm2}.
With $h=f_z$, $\mu_1= \mu/(1-|\nu|^2)$ and $\mu_2=\overline{\mu}\nu/(1-|\nu|^2)$, we may write this as 
$
h_{\overline{z}} = \mu_1 h_z +\mu_2 \overline{h_z} .
$
We have $|\mu_1|+|\mu_2| = |\mu|/(1-|\nu|) \leq || \mu/(1-|\nu|) ||_{L^{\infty}(\Omega)}$ almost everywhere. Write $k' = || \mu/(1-|\nu|) ||_{L^{\infty}(\Omega)}$. Since $|\mu|+|\nu|\leq k<1$ a.e., it follows that $k'\leq k$. 

Now \cite[Theorem 14.2.2, p.~369]{AIM}  implies that $h=f_z$ is $K'-$quasiregular, where $K'=(1+k')/(1-k')$, and then  \cite[Theorem 5.4.2, p.~175]{AIM} implies that $f_z$ belongs to  $W^{1,q}_{loc}(\Omega)$  for all $q\in (0,1+1/k')$. 
This proves part 4 of Theorem~\ref{thm2} and hence completes the proof of Theorem~\ref{thm2}. \hfill $\Box$

\bigskip

Of course the most natural way that the elliptic estimate (\ref{1}) is produced is from the nonlinear autonomous Beltrami equation $f_\zbar = {\cal A}(f_z)$ with the elliptic Lipschitz estimate 
\[  | {\cal A}(\zeta)-{\cal A}(\eta)|\leq k|\zeta-\eta|, \hskip10pt k<1,  \;\; \hbox{ for all }\zeta,\eta\in \IC   .  \]  
 Actually,  the method of ``frozen coefficients'' enables the nonlinear  Beltrami equation $f_\zbar = H(f,f_z)$ to be studied in this way as well,  see  \cite{AC2}.

\section{Proof of Theorem \ref{thm3}}

Let the assumptions of Theorem~\ref{thm3} be satisfied. 
It is clear each member of ${\cal F}$ is quasiregular with the same distortion bounds as for $f$.  We wish to prove a uniform bound for all $g\in {\mathcal F}$  in  $L^p(\ID)$, where $p>0$ is as in the assumptions of Theorem~\ref{thm3}.  We compute that
\begin{eqnarray*}
\int_\ID \Big| \frac{f(\lambda z)}{\lambda^{q/p}} \Big|^p \; \, dm(z) & = & \frac{1}{|\lambda|^{q}} \int_\ID | f(\lambda z)|^p \; \, dm(z) =  \frac{1}{ |\lambda|^{2+q}} \int_{\ID(0,|\lambda|)} | f(z)|^p \; \, dm(z)  \leq   C
\end{eqnarray*}
when $|\lambda|$ is small enough. So  we obtain a uniform bound in $L^p$ on the quasiregular family $\{ \frac{f(\lambda z)}{\lambda^s} : \lambda\in \ID\setminus\{0\} \}$. Now Theorem~\ref{Lpbounded} implies Theorem~\ref{thm3}.

\section{Quasiregular families bounded in $L^p$.}

We now present the proof of our main result,  Theorem \ref{Lpbounded}.

The conclusion of the theorem is already known to be valid if the elements of ${\mathcal F}$ are uniformly bounded in $L^{\infty}(\Omega)$. We also note that if $0<p<1$, then $L^p(\Omega)$ is not a Banach space,  but of course the $L^p$--norm $||f||_p = \left( \int_{\Omega} |f(z)|^p \,  dm(z)  \right)^{1/p}$ is well defined.

\medskip
 Let $\{f_k\}_{k=1}^{\infty}$ be a sequence in ${\mathcal F}$.  The Sto\"{\i}low factorisation theorem,  \cite[Theorem 5.5.1, p.~179]{AIM}, allows us to write
\begin{equation}
f_k = \varphi_k \circ g_k
\end{equation}
where each $g_k \colon \IC \to \IC$ is $K$--quasiconformal with $g_k(\infty) = \infty$ 
and $\varphi_k: g_k(\Omega)\to\IC$ is holomorphic.  Let $a_i$, $i=1,2$, be distinct points of $\Omega$ and let $\psi_k:\oC\to\oC$ be a M\"{o}bius transformation such that $\psi_k(g_k(a_i))=a_i$, $i=1,2$, and $\psi_k(\infty) = \infty$.  Writing
\[ f_k = \varphi_k\circ\psi_k^{-1} \circ \psi_{k} \circ g_k \]
we see that we can replace $g_k$ by $\psi_{k} \circ g_k$ so that we may and will in fact assume that each $g_k$ fixes two points of $\Omega$ and the point at infinity.  Then the family $\{g_k:\Omega\to\IC\}$ is precompact and passing to a subsequence without changing notation, we may assume that $g_k\to g$ locally uniformly in $\Omega$, where $g:\Omega\to \IC$ is $K$--quasiconformal (in fact, we may choose $g:\IC\to \IC$).  Similarly we write $\varphi_k$ for  $ \varphi_k\circ\psi_k^{-1}$.

In view of these properties of the functions $g_k$, to obtain the conclusion of the theorem it suffices to prove that the functions $\varphi_k$ form a normal family in every relatively compact subset of $g(\Omega)$, which we now proceed to do.

Let the domain $W$ be a relatively compact subset of $g(\Omega)$. Next, let the domain $U$ be a relatively compact subset of $\Omega$   such that $\overline{W} \subset g(U)$.  Then there is a domain $V\subset g(\Omega)$ such that $\overline{g(U)} \subset V$ and $\overline{V} $ is a compact subset of $g(\Omega)$. Since $g_k\to g$ uniformly on $U$, we have $W\subset g_k(U) \subset V$   for all large $k$; we may assume that this is true for all $k$.  

We denote the Jacobian determinant of $g$ by $J(w,g)$. Now each $g_k$ is $K$-quasiconformal in $\Omega$ and so by Astala's Theorem $J(w,g_k)\in L^q(U)$ for each $q \in [1,  \frac{K}{K-1})$.

We need the following lemma.

\begin{lemma}  \label{le1}
Let $\Omega$ be a planar domain and let $U$ be a relatively compact  domain in $\Omega$.  Suppose that for all $k\geq 1$, the mapping $g_k:\Omega \to\IC$ is $K$-quasiconformal, that $g_k:\Omega \to\IC$ is $K$-quasiconformal, and that  $g_k\to g$ locally uniformly on $\Omega$.   Then there is a constant $C=C(q,K)$ such that for every  $q$ with $1\leq q<\frac{K}{K-1}$, we have 
\begin{equation}
 \int_U   J(w, g_k)^q \; dm(w) \leq C  \int_U   J(w, g)^q \; dm(w).
\end{equation}
\end{lemma}

\noindent{\bf Proof.}  We suppose that $U=B$ is a relatively compact disk in $\Omega$,  the general result follows in an elementary manner.  Under the hypotheses we have $g_k\to g$ uniformly on $B$ and for all sufficiently large $k$,  $|g_k(B)|\leq |g(B)|+1<\infty$, where $|A|$ denotes the area of the set $A$.    The local uniform convergence on $\Omega$ implies weak convergence of the Jacobians,
\[ \int_\Omega \varphi(z) \; J(z,g_k) \, dm(z) \to  \int_\Omega \varphi(z) \; J(z,g) \, dm(z) \]
for every $\varphi\in C^{\infty}_{0}(\Omega)$. See
Corollary, p.~141, and Theorem 9.1, p.~216
in  \cite{S}.
We recall that if $J$ is the Jacobian of a $K$-quasiconformal mapping and $\omega= J^s$, where $\frac{-1}{K-1} <  s<\frac{K}{K-1}$,  then $\omega$ is an $A_p$-weight for all 
\begin{equation}
p> \left\{\begin{array}{ll} 1+s(K-1), & 0 \leq s <\frac{K}{K-1} \\ 1-\frac{s}{K} (K-1), & \frac{-1}{K-1} < s \leq 0 \end{array} \; \right.,
\end{equation}  
see \cite[Theorem 13.4.2]{AIM}.  With $s=q$ we have \cite[(13.56)]{AIM}
\begin{equation}
  \frac{1}{|B|} \, \int_B J(z,g_k)^q \, dm(z)  \leq  C(q,K) \left( \frac{|g_k(B)|}{|B|} \right)^{q} \leq C \left(\frac{|g(B)|+1}{|B|} \right)^{q} 
\end{equation}
where $C=C(q,K)$ is finite. Thus the sequence $\{J(z,g_k)\}$ is uniformly bounded in $L^q(B)$ and we may extract a weakly convergent subsequence in $L^q(B)$.  Since $J(z,g)$ also lies in $L^q(B)$ this weak limit must in fact be $J(z,g)$, and every weakly convergent subsequence will have the same limit.  The result follows.
 \hfill $\Box$
 
 \medskip

We can now use H\"{o}lder's inequality to see that if $q'= \frac{q}{q-1} > K$ and $s>0$, we have
\begin{eqnarray*}
 \int_{W} |\varphi_k(z)|^s\; \, dm(z)  & \leq &
 \int_{g_k(U)} |\varphi_k(z)|^s\; \, dm(z) 
 \\ & = &   \int_{g_k(U)} |(f_k \circ g_k^{-1} ) (z)|^s\; \, dm(z)
 \\  & = &  
  \int_U |f_k(w)|^s  J(w, g_k) \; dm(w)  
  \\ & \leq & \Big( \int_U |f_k(w)|^{sq'}  \; dm(w) \Big)^{1/q'} \; \Big( \int_U   J(w, g_k)^q \; dm(w) \Big)^{1/q} \\
 & \leq & C \Big( \int_U |f_k(w)|^{sq'}  \; dm(w) \Big)^{1/q'} \; \Big( \int_U   J(w, g)^q \; dm(w)  \Big)^{1/q} .  
\end{eqnarray*}
Thus if $f_k=\varphi_k\circ g_k$ is uniformly bounded in $L^{sq'}(U)$ for some $q'>K$, then the sequence $\{\varphi_k\}$ is uniformly bounded in $L^s(W)$. 

Now fix any $q$ with $1\leq q < K/(K-1)$. This determines $q'=\frac{q}{q-1} > K$. Then choose $s=p/q'>0$. Hence $sq'=p$. By the assumption of Theorem~\ref{Lpbounded}, we see that the numbers $\int_{W} |\varphi_k(z)|^s\; dz$ are uniformly bounded. Now it follows from Lemma~\ref{le2} below that the functions $\varphi_k(z)$ are locally uniformly bounded in $W$ and therefore that all the other claims of Theorem~\ref{Lpbounded} are valid.

\medskip

\begin{lemma}  \label{le2}
Let $\Omega$ be a domain in $\IC$. Suppose that $s$ and $M$ are positive real numbers. Let ${\mathcal F}$ be a family of holomorphic functions $f:\Omega\to\IC$ such that one of (a) and (b) below holds.

(a) For all $f\in {\mathcal F}$, we have
$$
\int_{\Omega} \log^+ |f(x+iy)| \, dx\, dy \leq M .
$$

(b) For all $f\in {\mathcal F}$, we have
$$
||f||_s = \left(  \int_{\Omega}  |f(x+iy)|^s \, dx\, dy \right)^{1/s} \leq M .
$$

Then the functions in ${\mathcal F}$ are locally uniformly bounded in $\Omega$, ${\mathcal F}$ is a normal family in $\Omega$, and the limit of every locally uniformly convergent sequence of functions in ${\mathcal F}$ is holomorphic in $\Omega$.
\end{lemma}

 \noindent{\bf Proof.}  Let $f\in {\mathcal F}$. 
 
Suppose that $b\in\Omega$ and set $d={\rm  dist} \, (b,\partial\Omega)>0$.
For $0<r<d$, write
$$
L(r,f) = \frac{1}{2\pi} \int_0^{2\pi} \log^+ |  f(b+re^{i\theta})  |  \, d\theta ,
$$
$$
I_s(r,f) = \frac{1}{2\pi}  \int_0^{2\pi} |  f(b+re^{i\theta})  |^s  \, d\theta .
$$ 
Since $\log^+ |f|$ and $|f|^s$ are subharmonic in the disk ${\mathbb D} (b,d) = \{z: |z-b|<d \}$, we have
$$
\log^+ |f(b) | \leq L(r,f) ,
$$
$$
| f(b) |^s \leq I_s(r,f)
$$
for every $r\in (0,d)$. 
Also, since
$\log^+ |f|$ and $|f|^s$ are subharmonic in  ${\mathbb D} (b,d) $, the functions $L(r,f)$ and $I_s(r,f)$ are increasing functions of $r$ for $0<r<d$. Denote the annulus $\{z: d/2<|z-b|<d \}  $ by $A$. Applying each of the above inequalities with $r=d/2$, and then multiplying both sides by $r$ and integrating with respect to $r$ from $d/2$ to $d$ we get
\begin{eqnarray*}
(3/8)  d^2 \log^+ |f(b) |  &=& (1/2) (d^2-(d/2)^2) \log^+ |f(b) |  \\  & \leq &  
\int_{d/2}^d  L(d/2,f) \, r\, dr \leq \int_{d/2}^d  L(r,f) \, r\, dr
\\  &=&
 \frac{1}{2\pi}  \int_A \log^+ |f(x+iy)| \, dx\, dy \leq \frac{M}{2\pi}  
\end{eqnarray*}
if the assumption (a) is satisfied, and similarly 
\begin{eqnarray*}
(3/8)  d^2 | f(b) |^s & \leq &  \int_{d/2}^d  I_s(r,f) \, r\, dr
\\  &=&
 \frac{1}{2\pi}  \int_A  |f(x+iy)|^s \, dx\, dy \leq \frac{M^s} {2\pi}  
\end{eqnarray*}
if the assumption (b) is satisfied. 

Thus
$$
|f(b) | \leq \exp \left( \frac { 4 M } { 3 \pi d^2  }  \right)
$$
in case (a), and
$$
|f(b) | \leq M \left( \frac { 4 } { 3 \pi d^2  }  \right)^{1/s}
$$
in case (b). 
This proves that the functions in ${\mathcal F}$ are locally uniformly bounded in $\Omega$, and the remaining claims now follow from standard results in complex analysis. This completes the proof of Lemma~\ref{le2}. \hfill $\Box$

\bigskip
 \noindent A. Hinkkanen, University of Illinois at Urbana-Champaign,   aimo@math.uiuc.edu\\ 
 \noindent G. Martin, Massey University,   New Zealand,  g.j.martin@massey.ac.nz

\end{document}